\documentclass[11pt]{article}
\pdfoutput=1
\usepackage{epsfig}

\usepackage{enumerate}

\begin{document}
\newcommand{\arcsech}{\mbox{arcsech}}
\newcommand{\bc}{\begin{center}}
\newcommand{\ec}{\end{center}}
\newcommand{\be}{\begin{equation}}
\newcommand{\ee}{\end{equation}}
\newcommand{\bea}{\begin{eqnarray}}
\newcommand{\eea}{\end{eqnarray}}
\newcommand{\bean}{\begin{eqnarray*}}
\newcommand{\eean}{\end{eqnarray*}}
\newcommand{\bd}{\begin{description}}
\newcommand{\ed}{\end{description}}

\newcommand{\caln}{\cal N}
\newcommand{\caloi}{{\cal OI}}
\newcommand{\hatv}{{\bf v}}
\newcommand{\hatrho}{\hat{{\bf \rho}}}
\newcommand{\hatlam}{\lambda}
\newcommand{\tillam}{\tilde{\lambda}}
\newcommand{\tT}{\tilde{T}}
\newcommand{\hT}{\hat{T}}
\newcommand{\bT}{{\bf T}}
\newcommand{\bA}{{\bf A}}
\newcommand{\ep}{\epsilon}
\newcommand{\cd}{{\cal D}}
\newcommand{\hp}{\hat{T}}
\newcommand{\epf}{\epsilon_f}
\newcommand{\Nn}{N_{net}}
\newcommand{\Ns}{N_{sph}}
\newcommand{\Nd}{N_{dat}}
\newcommand{\lng}{\langle}
\newcommand{\rng}{\rangle}
\newcommand{\vep}{\epsilon}
\newcommand{\vph}{\varphi}
\newcommand{\vp}{\varphi}
\newcommand{\bfk}{{\bf k}}
\newcommand{\bfr}{{\bf r}}
\newcommand{\bfN}{{\bf N}}
\newcommand{\bfa}{{\bf a}}
\newcommand{\bfb}{{\bf b}}
\newcommand{\bfc}{{\bf c}}
\newcommand{\bfH}{{\bf H}}
\newcommand{\bfK}{{\bf K}}
\newcommand{\bfL}{{\bf L}}
\newcommand{\bfM}{{\bf M}}
\newcommand{\bfx}{{\bf x}}
\newcommand{\bC}{{\bf C}}
\newcommand{\bp}{{\bf{p}}}
\newcommand{\bfe}{{\bf e}}
\newcommand{\hbt}{\hat{\bar{T}}}
\newcommand{\uvr}{{\bf \hat{r}}}
\newcommand{\uvn}{{\bf \hat{n}}}
\newcommand{\bw}{\bar{w}}
\newcommand{\bx}{\mbox{\boldmath $\xi$}}
\newcommand{\bnu}{\mbox{\boldmath $\nu$}}
\newcommand{\my}{M_{yr}}

\newcommand{\bvt}{\begin{verbatim}}
\newcommand{\evt}{\end{verbatim}}
\newcommand{\dbar}{\vrule width 0.25truein height 1.0pt depth -0.6pt}

\renewcommand{\baselinestretch}{2.0}

\vskip 0.5cm

\begin{center}
{\bf\Huge DD Calculus}\\
\vskip 0.5cm
{\huge Learning Direct Calculus in Two Hours\footnote[1]{This is an English translation with comments of a paper in Chinese:	Shen, S.S.P., and Q. Lin, 2014: Two hours of simplified teaching materials for direct calculus, {\it Mathematics Teaching and Learning}, {\bf No. 2}, 2-1 - 2-6.   The Chinese original is attached at the end of this translation. 
}}
\end{center}
\vskip 1.5cm

\renewcommand{\baselinestretch}{1.2}

\begin{center}
{\Large SAMUEL S.P. SHEN }  \\
Department of Mathematics and Statistics, 
San Diego State University,
San Diego, CA 92182, USA. Email: 
{\tt sam.shen@sdsu.edu}\\
\vskip 0.5cm
{\Large QUN LIN}\\
Institute of Computational Mathematics,
Chinese Academy of Sciences,
Beijing 100080, CHINA.  Email: 
{\tt linq@lsec.cc.ac.cn}\\
\vskip 0.5cm
 
\end{center}

\thispagestyle{empty}
\renewcommand{\baselinestretch}{1.2}

\vskip 0.9cm
\noindent Citation of this work: Shen, S.S.P., and Q. Lin, 2014: Two hours of simplified teaching materials for direct calculus, {\it Mathematics Teaching and Learning}, {\bf No. 2}, 2-1 - 2-6.
\vskip 0.9cm

\noindent {{\bf  Summary}} This paper  introduces DD calculus and describes the basic  calculus concepts of derivative and integral in a 
direct and non-traditional way, without limit definition: Derivative is computed from the point-slope 
equation of a tangent line and integral is defined as the height increment of a curve. This direct approach to calculus has three distinct features: (i) it defines derivative and 
(definite) integral without using limits, (ii) it defines derivative and antiderivative  simultaneously via 
a derivative-antiderivative (DA) pair, and (iii) it posits the fundamental theorem of 
calculus as a natural corollary of the definitions of derivative and integral. The first D in DD calculus attributes to Descartes for his method of tangents and the second 
D to DA-pair. The  DD calculus, or simply direct calculus, makes many traditional notations and procedures 
unnecessary, a plus when introducing calculus to the non-mathematics majors. It has few intermediate procedures, which can help 
dispel the mystery of calculus as perceived by the general public. The materials in  this paper are intended for use in a two-hour introductory lecture on calculus.

\newpage

\section{Introduction: Necessity to dispel calculus mystery and simplify calculus notations}

\noindent  ``Calculus" is not a commonly used word in daily life. The Oxford Dictionary indicates that the word comes from the mid-17th century  Latin and literally 
means small pebble (such as those used on an abacus) for counting. The dictionary gives three meanings of ``calculus": a branch of mathematics that deals with derivatives
and integrals, 
a particular method of calculation, and a hard mass formed by minerals. Obviously  here we are interested in the first meaning. 
We wish to demonstrate that (i) the calculus method can be developed by analyzing steepness and height change of a curve,
(ii) the method development can be achieved directly using Descartes' method of tangents and does not need an introduction of limit as prerequisite,  
and (iii) the basic method of calculus and a few simple examples
can be introduced in a one-hour or two-hour lecture to a high-school level audience. 

 Calculus is one of the most important tools in a knowledge based society. 
Millions of people around the world learn calculus everyday. All engineering, science, and business major undergraduate students 
must  take calculus. Many high schools offer calculus courses. The usefulness and power of calculus have been well recognized. 
Nonetheless,  calculus is  a mysterious subject to many people and is regarded by the general public  as accessible only to a few privileged 
people with special talents. Tight schedules and high fail rates for the first semester calculus have given the course 
a reputation as a monster, a nightmare, or a phycological barrier for many students, some of whom are even  STEM (science,
technology, engineering and mathematics) majors.  Calculus can be a topic that causes people at a social gathering to shake their heads in incomprehension, 
shy away from the daunting challenge of understanding it, or express effusive exclamation of awe and admiration. 
It is also sometimes associated with conspicuous nerdiness. In classrooms, the student-instructor 
relationship can be tense. Some students regard  calculus instructors as 
inhuman and ruthless aliens,  while instructors frequently joke about students' stupidity, clumsiness, or silly errors.  Tedious and peculiar notations coupled
 with  fiendish and complex approaches to calculus teaching and learning
  may have contributed to the above unfortunate situation. 

A major cause of this mystery and scare of calculus is unnecessarily complex terminologies, including definite and indefinite integrals,
derivatives defined as limits, definite integrals defined as limits, the difference between $dx$ and $\Delta x$, and using ever-finer divisions 
of an area under a curve to approximate a definite integral (i.e., introducing the definite integral by area and limit), to name but a few.
Additionally,  the Riemann sum 
with its arbitrary point $x_i^*$ in the interval $[x_i, x_{i+1}]$ complicates calculation procedures and  adds more confusion. These conventional concepts 
and notations are essential for professional mathematicians who research mathematical analysis, but are absolutely unnecessary to 
the majority of calculus learners, who are majoring in engineering, science, business, or other non-mathematical or non-statistical fields. 

In the first semester calculus, instructors repeatedly emphasize that an indefinite integral yields a function while a definite integral yields
a real value. When an instructor discovers that some students still cannot tell the difference between a definite integral and an indefinite integral in the final exam,
s/he becomes disappointed and complains that these terrible students did not pay attention to her/his repeated emphasis on the difference,
but s/he rarely questions the necessity of introducing the concept of indefinite integrals and their notations.

The purpose of this paper is to dispel this mystery of calculus   by introducing a more direct approach. 
We attempt to introduce  the basic ideas of calculus in one hour without using the concept of limits. 
Our introduction to derivatives directly uses the idea of Ren{\' e} Descartes' (1596-1650) method of tangents (Cajori,1985, pp.176-177;  Coolidge, 1951; Susuki, 2005; Range, 2011), rather than indirectly uses the 
 secant line method with a limit. 
We introduce derivatives and antiderivatives simultaneously, using derivative-antiderivative (DA) pairs. 
Our introduction to integrals is directly from the DA pairs and the height increment of an antiderivative curve. The height increment approach 
has been advocated by Q. Lin in China for over two decades  (see Lin (2010) and Lin (2009) for two recent examples).  We 
define the area under the curve of an integrand by the integral, and then explain why the definition is reasonable. This is 
a reversal of the traditional definition, which defines an integral by  calculating the area underneath the  curve of an  integrand.

We also demonstrate that an introduction to the basic ideas of calculus does not need to use many complicated notations. Thus, the notations of 
derivatives $f'(x)$ and integrals $I[f(x), a,b]$  in this paper are parsimonious,  simple, computer friendly and non-traditional. 

In the following sections, we will first introduce (a) the method of point-slope equation for 
calculating derivatives, and (b)  DA pairs for calculating both antiderivatives and derivatives. Then we will  introduce integrals as 
the height increment of  an antiderivative curve. Next we will discuss the mathematical rigor of our direct calculus
method, but this part of material does not need to be included in a one-hour lecture. Finally we provide a brief history on the 
development of calculus ideas and give our conclusions on the direct calculus from the perspective of Descartes' method of tangents and DA pairs.

\section{Slope, and derivative and antiderivative pairs }

\noindent When we drive on a steep highway, we often see a grade warning sign like the one in Figure 1. The 9$\%$ grade or slope  on a highway means 
that the elevation will decrease 90 feet when the horizontal 
distance increases 1000 ft. The grade or slope is calculated by
\be
m=\tan \theta = \frac{H}{L}=\frac{90}{1000}
\ee
i.e., the ratio of the opposite side to the adjacent side of the right triangle in Figure 1.

 \begin{figure} [ht]
\centering
\includegraphics[height=2.6in, width=6in]{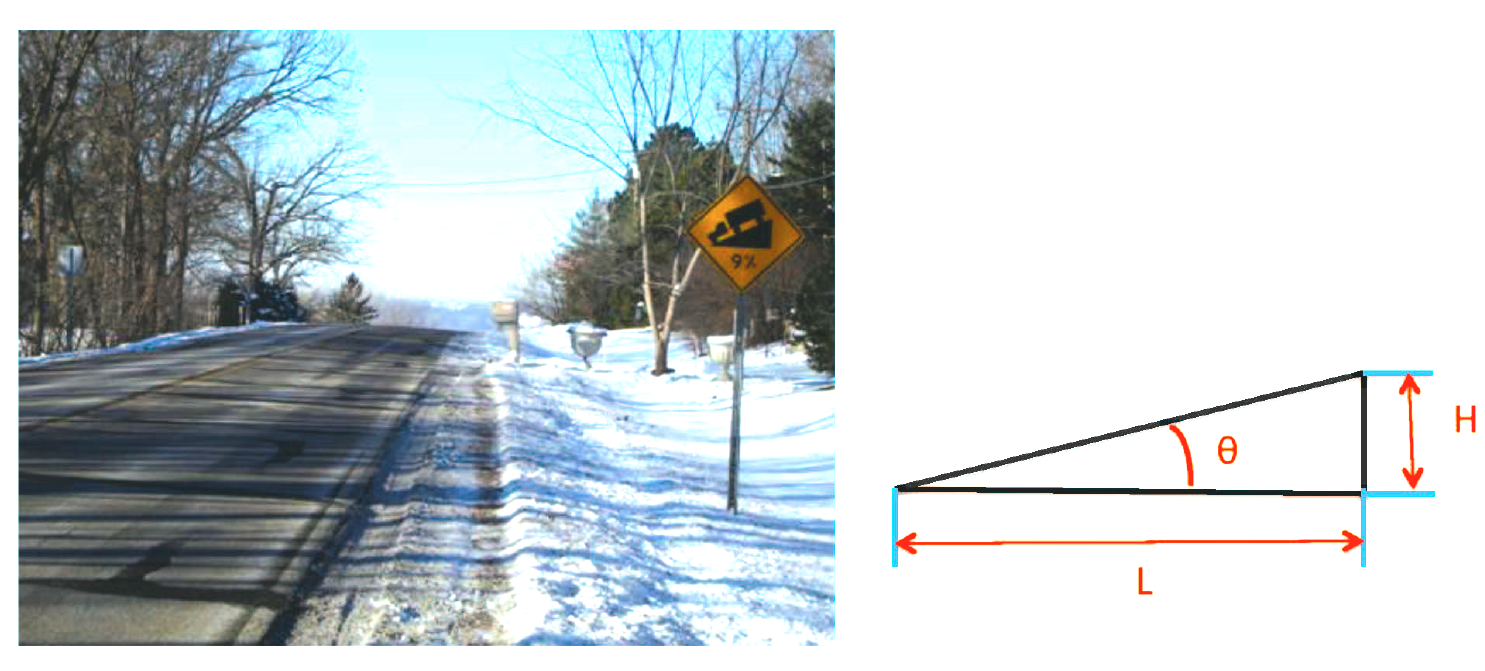}
\caption{ Highway grade sign and a right triangle to show slope.}
\label{figure1}
\end{figure}

The slope of a curve at a given point is defined as the slope of a tangent line at this point. 
Figure 2 shows three points: P, A and B. 
$T_P$ represents the tangent line at $P$ whose slope is defined as the slope of the curve at P. 
The  tangent line's slope is used to measure
 the curve's steepness. Calculus studies (i) the slope $\tan \theta$  of a curve at various points, and (ii) the height 
increment $H$  from one point to another, say, A to B, as shown in Figure 2. Our geometric intuition 
indicates that height $H$ and slope $\tan \theta$ are related because  $H$ increases
rapidly  if the slope is large for an upward trend. A core formula of calculus is to describe the relationship between
these two quantities.

\begin{figure} [ht]
\centering
\includegraphics[height=2.5in, width=3in]{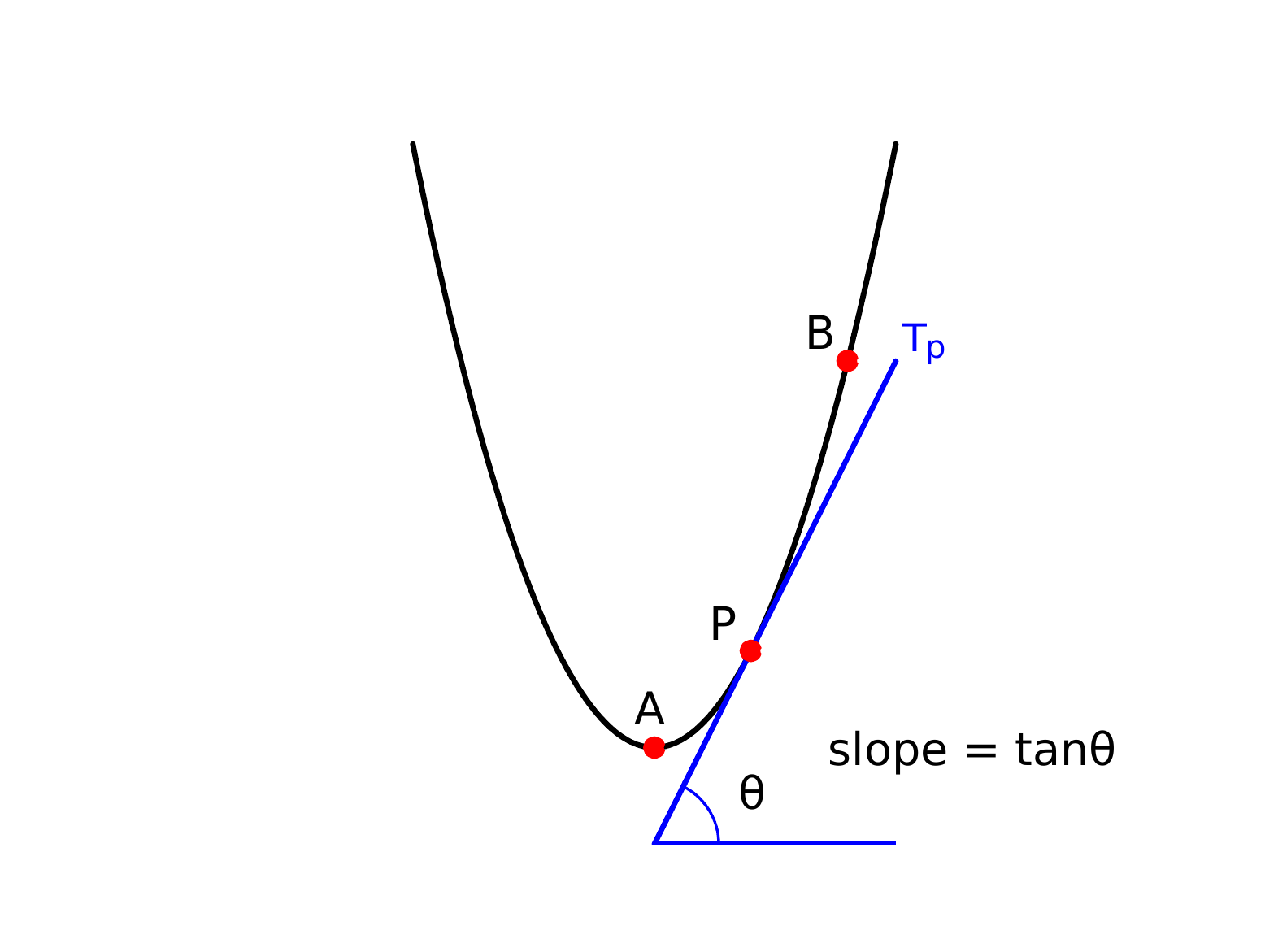}
\caption{ A curve, three points and a tangent line without coordinates.}
\label{figure2}
\end{figure}

Let us introduce the coordinates x and y and use a function $y=f(x)$ to describe the curve. We start with an example
$f(x) = x^2$.

The tangent line of the curve at point $P(x_0, y_0)$ can be described by a point-slope equation
\be
y-y_0=m(x-x_0), \label{eq2}
\ee
where $y_0=x_0^2$, $m$ is the slope to be determined by the condition that the tangent line touches the curve at 
one point. In fact, ``tangent" is a word derived  from Latin and means ``touch." 

The tangent line (\ref{eq2}) and the curve
\be
y=x^2 \label{eq3}
\ee
have a common point $P(x_0, y_0)$ (See Figure 3), where $x_0$ must be a double root, since 
the straight line is a tangent line.

\begin{figure} [ht]
\centering
\includegraphics[height=2.8in, width=3in]{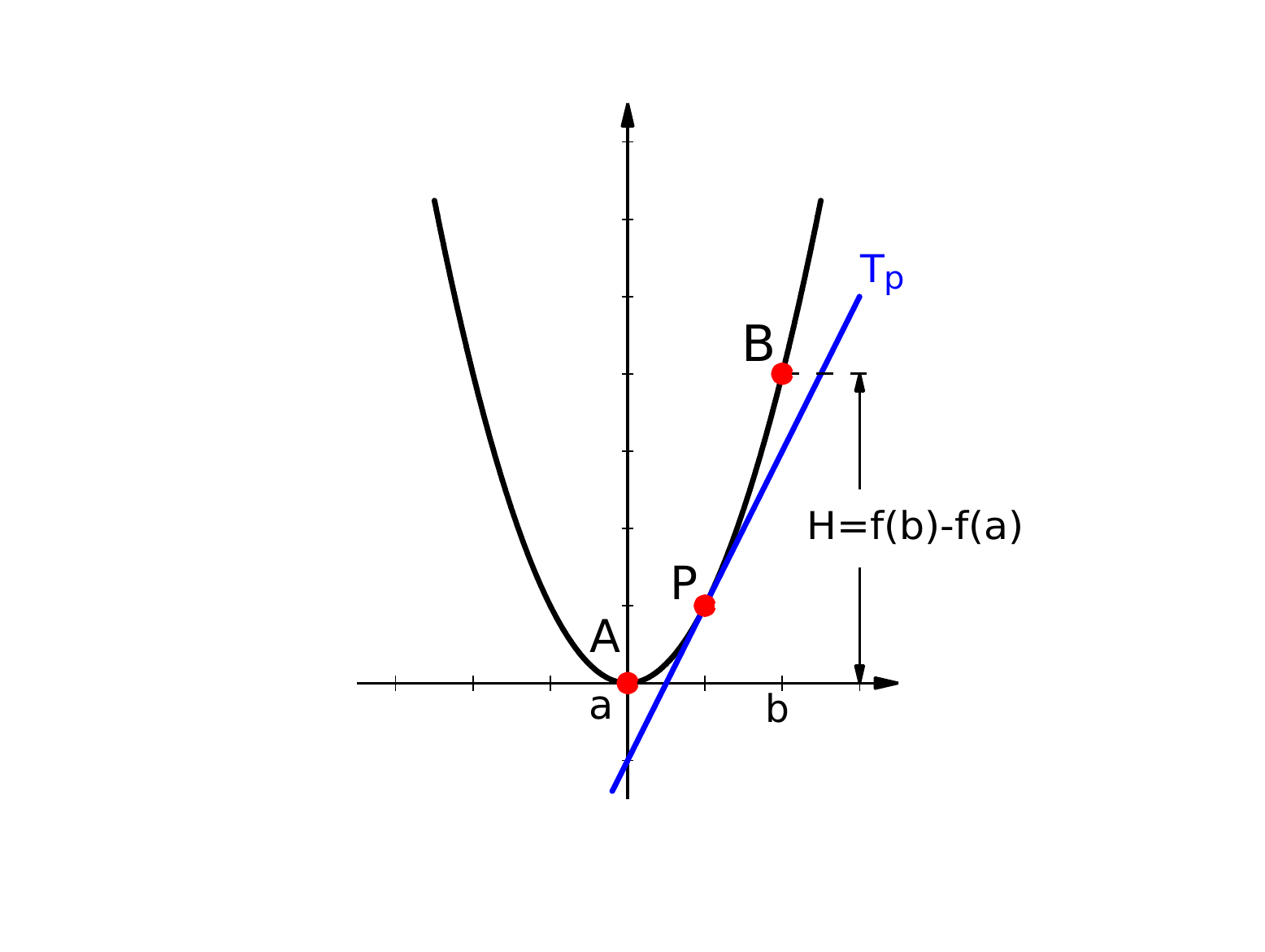}
\caption{A curve, three points and a tangent line on xy-plane.}
\label{figure3}
\end{figure}

Substituting (\ref{eq3}) to (\ref{eq2}) to eliminate $y$, we have 
\be
x^2 - x_0^2 = m(x-x_0),
\ee
then
\be
(x-x_0)(x+x_0)-m(x-x_0)=0,
\ee
or 
\be
(x-x_0)(x+x_0-m)=0.
\ee
The two solutions of this quadratic equation are
\be 
x_1 -x_0=0
\ee
and 
\be
x_2+x_0 -m=0
\ee
Since $P(x_0,y_0)$ is a tangent point, $x_0$ must be a repeated root, also called a double root. Hence, 
$x_1=x_2=x_0$. This yields 
\be
m=2x_0.
\ee

Thus, we claim that the slope of the curve $y=x^2$ at $x_0$ is $2x_0$, at $a$ is $2a$, at $3$ is $2\times 3=6$, and in general, at $x$ is $2x$. 
The slope measures the steepness of the curve $y=f(x)$, i.e., the rate of  the curve's height increase or decrease. 
The slope, or rate, varies from point to point. The slope is thus a derived quantity from the original function $y=f(x)$ and 
is called ``derivative." We have the following definition.

\vskip 0.2in

\noindent {\bf Definition 1. (Definitions of derivative as slope and DA pair).}
{\it The  slope  $2x$ is called the derivative of $x^2$. Further, $x^2$ is called an antiderivative of $2x$. And ($2x$, $x^2$) is called
a derivative-antiderivative (DA) pair. }

\vskip 0.2in

For a general function $y=f(x)$,  the slope of the  curve $y=f(x)$ is called {\it derivative} and is denoted by $f'(x)$ or $y'$, and $f(x)$ is 
called antiderivative of $f'(x)$. Thus, $(f'(x), f(x))$ is a DA pair. 

If $f(x)=C$ is a constant, then $y=C$ represents a horizontal line whose slope is $0$ for any $x$. Hence $(C)'=0$, and $(0, C)$ 
is a DA pair. 

 If $f(x)$ is a linear function,  
then $y = \alpha + \beta x$ represents a straight line whose slope is $\beta$ for any $x$, 
hence $(\alpha + \beta x)'=\beta$, i.e., $(\beta, \alpha + \beta x)$ is a DA pair. 

Thus, the derivative's geometric meaning is the slope of the curve $y=f(x)$: $f'(x)$ is large at places where the curve $y=f(x)$ is steep. 
At a flat point, such as the maximum or minimum point of $f(x)$, the slope is zero since the tangent lines at these 
points are horizontal. 

In addition to the geometric meaning, derivative has physical meaning, such as  speed, and biological meaning, such as  growth rate, as well as 
the meaning of the rate of change in almost any scientific  field and anyone's daily life. 
As an example, for a car driven at $v=60$ mph (miles per hour) for two hours, the total distance 
traveled is  $s=v\times t= 60\times 2 = 120$ mi. Here, $(v, vt)$ or $(v, s)$ is a DA pair for a general time $t$. 

Free fall is another example.  An object's free fall has its distance of falling equal to $s=(1/2)gt^2$ and its falling 
speed is $v=gt$, where $g=32 [ft/s^2]$ is the Earth's gravitational 
acceleration. Galileo Galilei (1564-1642) discovered this time-square relationship for the distance.
Since derivative $t^2$ with respect to $t$ is $2t$, we have $s'=\frac{1}{2}g (2t)=gt$. Thus, $(gt, (1/2)gt^2)$ or $(v,s)$ is a DA pair. 

In general, the meaning of derivative is the rate of change of the function $f(x)$ with respect to the independent variable $x$, which 
can be either time or spatial location. 

Let us now return  to the derivative calculation. The above tangent line approach for finding the slope for $x^2$ can be applied to the function $y=x^3$. It is to solve the following simultaneous equations 
\bea 
&& y-y_0 = m (x-x_0),\\
&&  y=x^3.
\eea
Eliminating $y$, we have 
\be
x^3-x_0^3 = m(x-x_0).
\ee
The factorization of this equation yields
\be
(x-x_0)(x^2+xx_0 + x_0^2 -m)=0.
\ee
This factorization implies $x_1-x_0=0$ and $x_2^2 +x_2x_0+x_0^2 -m=0$. 
The cubic equation has three solutions. Because $P(x_0,y_0)$ is the  tangent point, $x_0$ is a repeated root of these two equations: $x_1=x_2=x_0$,
which leads to 
\be
m=3x_0^2.
\ee
Hence, we claim that the derivative of $x^3$ is $3x^2$, and $x^3$ is an antiderivative of $3x^2$. Namely, ($3x^2$,  $x^3$) forms a 
DA pair. 

Following the tangent line approach, the above  examples have demonstrated four  DA pairs:
\be
(0,C),\\
(1,x),\\
(2x, x^2), \\
(3x^2, x^3).
\ee
The tangent line approach can be applied to any power function $x^n$, where $n$ is a positive integer. The DA pair for $x^n$ is 
\be
(nx^{n-1}, x^n).
\ee
This formula actually holds for any real number $n$ with the exception of $n=0$. For example, $(x^{1/2})'=(1/2)x^{-1/2}$. 
Proof of this claim is not easy, but fortunately the calculation of derivatives and antiderivatives can be easily 
done using computer programs. Many open source computer software packages are available to do this kind of calculation. 
For example, WolframAplha is one. At the website 
www.wolframalpha.com, 
you can enter a derivative command and a function, such as 
\begin{verbatim}
derivative  x^(3/2)+4x^3
\end{verbatim}
The computer will give you its derivative, plus other information, such as the graph of the 
derivative function  (see Figure 4). You can also use WolframAlpha via a smart phone application.

\begin{figure} [ht]
\centering
\includegraphics[height=4in,width=5in]{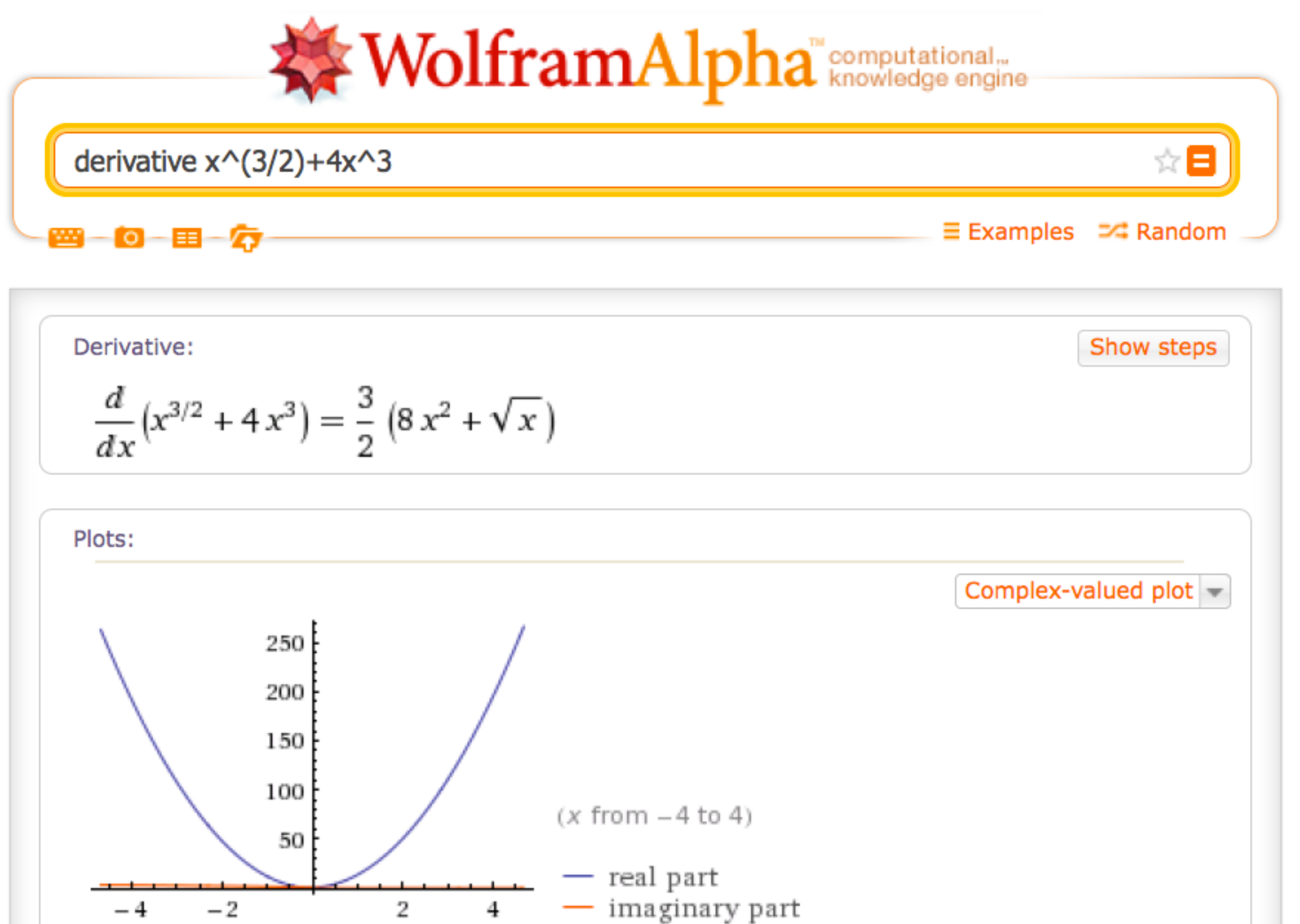}
\caption{WalframAlpha derivative example.}
\label{figure4}
\end{figure}

To find an antiderivative, you can use a similar command
\begin{verbatim}
antiderivative  x^(3/2)+4x^3
\end{verbatim}

Using this program, one can easily find DA pairs for commonly used functions. See the list below.
\begin{enumerate}[(i)]
\item Exponential function: $(e^x, e^x)$.
\item Natural logarithmic function: $(\frac{1}{x}, \ln x)$.
\item Sine function: $(\cos x, \sin x)$.
\item Cosine function: $(-\sin x, \cos x)$.
\item Tangent function: $(\sec^2 x, \tan x)$.
\end{enumerate}

\section{ Height increment and integrals}

 \noindent When we trace a curve, we care about not only the slope, but also the ups and downs of the curve, i.e., the 
increment or decrement of the curve from one point to another. When we drive over a mountain road,
we also care about both steepness (i.e., slope)  and elevation. Apparently, the slope and height increment are related. The slope has 
already been defined as derivative in the above section. In this section, the height increment is defined as {\it integral}, 
since the height increment or elevation increment is an integration process, or an accumulation process of point 
motion, measured by both speed and time. 

For a function $y=f(x)$,  its increment from $A (a, f(a))$ to $B (b, f(b))$ is $f(b)-f(a)$ as shown in Figure 3. Another notation for the increment is $f(b)-f(a) =f(x)|^b_a$.
This height increment is used to describe the  integral 
definition below. 

\vskip 0.2in 

\noindent {\bf Definition 2. (Definition of integral as height increment of a curve). }{\it The function's increment $f(b)-f(a)$ from $A (a, f(a))$ to $B (b, f(b))$ is  defined as the integral of the derivative function $f'(x)$ in the interval $[a,b]$ and is denoted by $I[f'(x), a, b]=f(b)-f(a) $. Here, $f'(x)$ is called the integrand, and $[a,b]$ is called 
the integration interval. }

\vskip 0.3cm

\noindent {\bf Example 1.}
Given $f(x)=x$,  $f'(x)=1$, and $[a,b]=[0,2]$, we have
\be
I[f'(x), a, b]=I[1, 0,2] =x|^2_0= 2-0 = 2.
\ee
The area between $y=1$ and $y=0$ in the interval $[0,2]$ is also 2 (see Figure 5 for $f'(x)=1, a=0$, and $b=1$).

\begin{figure} [ht]
\centering
\includegraphics[height=2.2in,width=4in]{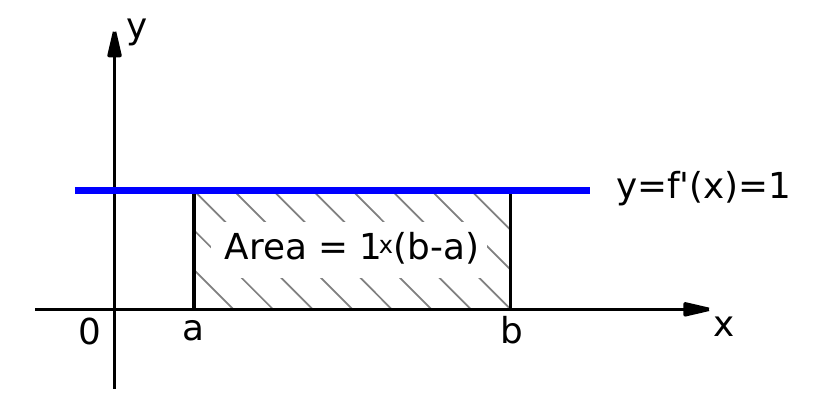}
\caption{The area of a rectangle under a horizontal line.}
\label{figure5}
\end{figure}

\vskip 0.3cm

\noindent {\bf Example 2. } If we integrate speed $v(t)$, we will then get the distance $I[v(t), a, b]$ travelled from time $t=a$ to $t=b$. 
If $v(t)$ is a constant, say, $60$ mi/hour, and if $a=14:00$ and $b=16:00$, then the integral $I[60, 14,16]=60t|_{14}^{16}=60\times (16-14)= 120$ mi is the 
total distance traveled from 2pm to 4pm. Here $60t$ is an antiderivative of $60$.  If we plot $v$ as a function of $t$, then $120$ is equal to the area 
of the rectangle bounded by $v=60, v=0, t=14$ and $t=16$ (See Figure 5 for $f'(x)=60$, $a=14$, and $b=16$). 

\vskip 0.3cm

\noindent {\bf Example 3. }
 Given $f(x)=(1/2)x^2$,  $f'(x)=x$, and $[a,b]=[0,1]$, we have
\be
I[f'(x), a, b]=I[x, 0,1] = (1/2)x^2|_0^1=(1/2)(1^2-0^2 )= 1/2.
\ee

The area under the integrand $y=x$ but above the $x$-axis in $[0,1]$ is $1/2$ (see Figure 6 for $a=0$ and $b=1$).

\begin{figure} [ht]
\centering
\includegraphics[height=3in,width=3in]{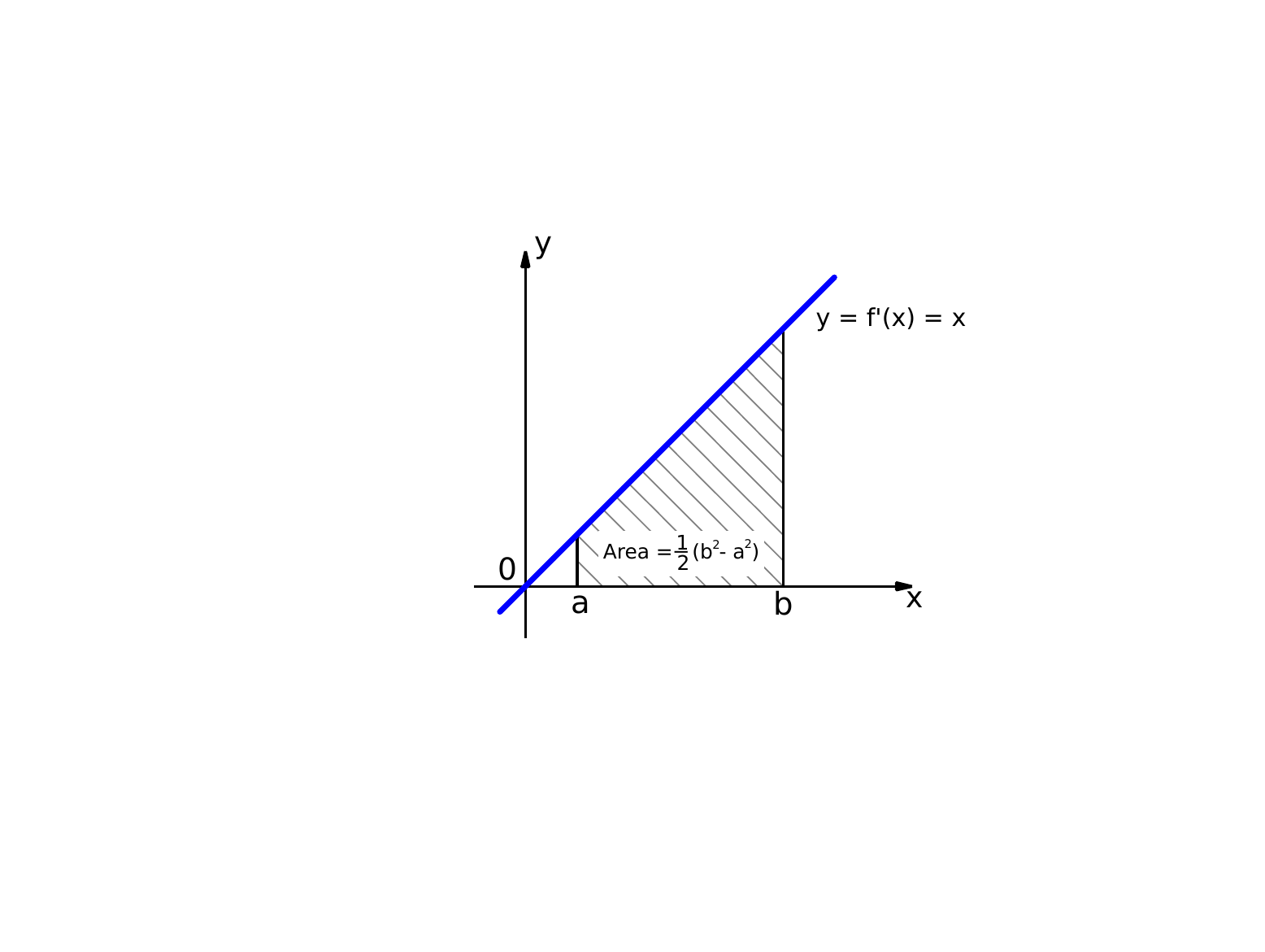}
\caption{The area of a triangle under a straight line. }
\label{figure6}
\end{figure}

\vskip 0.3cm

\noindent {\bf Example 4. } In the free fall problem, the speed is a linear function of time, $v=gt$, and the 
integration $I[gt,0,x]= (1/2)gx^2$ is the distance traveled from time zero to time $x$. The region bounded 
by $v=gt, v=0, t=0$ and $t=x$ is a right triangle with base equal to $x$, height $gx$ and area $(1/2)\times x\times gx=(1/2)gx^2$.  In 
this example, we have chosen to use $x$ as an arbitrary right bound for the region. This $x$ can can be any number, such as 1, 2, or 2.5. 

\vskip 0.3cm

In the above four examples, the area under a curve is equal to an integral. 
As a matter of fact, this inference of area equal to an integral is generally true. For an irregular region, we can simply use 
the integral $I[f'(x),a,b]$ as the definition of the area of the region bounded by $y=f'(x)$,  $x$-axis,  $x=a$ and $x=b$. 
The next section will justify this definition. As for the calculation of an integral, if one knows the relevant  DA pair,
the integral is a simple substitution $f(b)-f(a)$. Otherwise, one can use a computer to find the antiderivative, or to directly  evaluate
$I[f'(x), a,b]$. 
 Similar to derivative software, there are many free  online computer programs and smart phone apps  for calculating 
integrals. Figure 7 shows the 
WolframAlpha calculation of the integral $I[x^4+2x,0,1]$ using the  command  
\begin{verbatim}integrate x^4+2x from 0 to 1
\end{verbatim}
The result is 6/5. i.e., $I[x^4+2x,0,1]=6/5$.

\begin{figure} [ht]
\centering
\includegraphics[height=4in,width=5.5in]{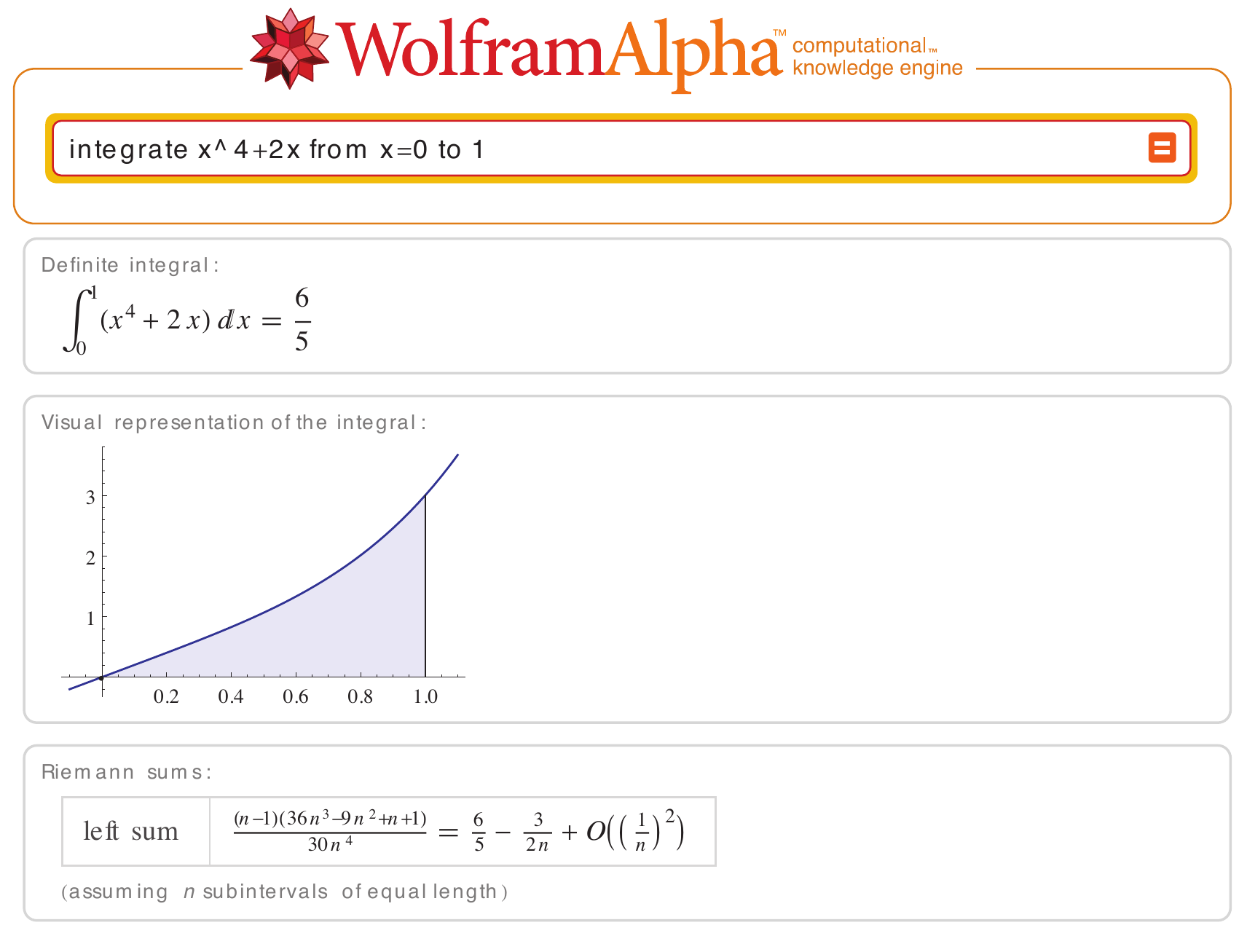}
\caption{WalframAlpha integral example. }
\label{figure7}
\end{figure}

The definition of an integral states that the integral of a function in an interval is the increment of its antiderivative in the same interval. One can 
write 
\be
I[g(t), a,b] = G(b)-G(a)
\ee
where $G(t)$ is an antiderivative of $g(t)$. Another way to express the above is 
\be
I[G'(u), a,b] = G(b)-G(a).
\ee
In the above two expressions, $t$ and $u$ are the integration variables, also called dummy variables.
The integral values are independent of the choice of dummy variables.  One can use any symbol to 
represent this variable. In practical applications, if the independent variable is time, such as when speed is a function of time, 
$t$ is often used as the independent variable. 

Also according to the integral definition, the integral of $f'(t)$ in the interval $[a,x]$ is 
\be
I[f'(t),a,x] = f(x)-f(a).
\ee
Taking derivative of both sides of this equation with respect to $x$, we have 
\be
(I[f'(t),a,x] )'_x= f'(x)  \label{ftc2}
\ee
since $(f(a))'_x=0$ due to  $f(a)$ being a constant with respect to $x$ and  having  a slope of zero. Here we have used a subscript $x$ to indicate 
that the independent variable is $x$ and the derivative is with respect to $x$. 

Equation (\ref{ftc2}) is often called Part II of the Fundamental Theorem of Calculus (FTC), while 
the definition of an integral is actually often referred to as Part I of the FTC. Part II of the FTC is saying that an antiderivative 
can be  explicitly expressed  by an integral. Thus, the FTC makes a computable and close connection between slope and height increment,
and enhances our intuitive sense that the height increment in an interval is closely  related to the slope of our interested curve, i.e., 
our study function. 

\vskip 0.3cm

\noindent {\bf Example 5. } $I[x,0,1] = (x^2/2)|^1_0 = 1^2/2-0^2/2=1/2$, because  $(x,x^2/2)$ is a DA pair. The area between $y=x$ and $y=0$ over 
$[0,1]$ is 1/2 (=$I[x, 0,1]$) (See Figure 6).

\vskip 0.3cm

\noindent {\bf Example 6.} $I[x^2,0,1]= (x^3/3)|_0^1=1^3/3-0^3/3=1/3$, because $(x^2, x^3/3) $ is a DA pair, or simply since 
$x^2$'s antiderivative is $x^3/3$.  The area of the right triangle with a curved hypotenuse bound by $y=x^2, y=0$ and $x=1$ is 
$1/3$ (=$I[x^2, 0, 1]$) (See Figure 8). The WolframAlpha command for this calculation  is {\tt integrate x}\^{\tt 3 from 0 to 1}. 

\begin{figure} [ht]
\centering
\includegraphics[height=3in,width=3in]{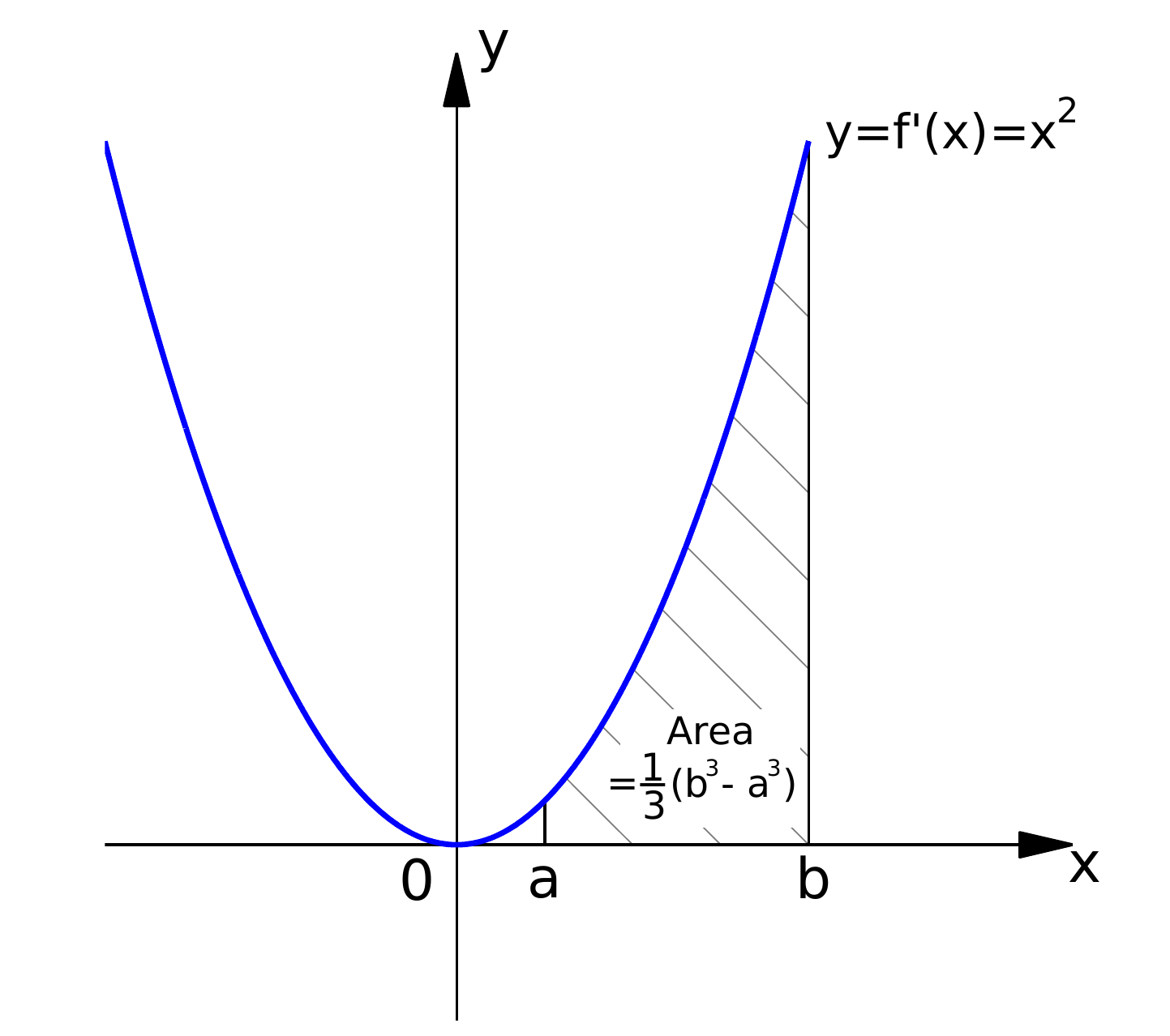}
\caption{ The area of a triangle of a curved hypotenuse when $a=0$.  }
\label{figure8}
\end{figure}

\section{ Discussion and mathematical rigor of direct calculus}

\noindent Two points  are discussed here. First, is our definition of area by an integral reasonable and mathematically rigorous? Second, 
besides using computer programs to calculate the DA pairs of complicated functions, can one 
provide a systematic procedure of hand calculation? 

First, how do we know that our  definition of area using an integral is reasonable? According to the Oxford dictionary,``area" is defined 
as ``the extent or measurement of a surface or piece of land." The word ``area" comes from the mid-16th century  Latin, literally meaning a ``vacant piece of 
level ground."  The units of an area are ``square feet",``square meters", etc, meaning that the area  of a region is equal to the number of equivalent squares,  
each side equal to a foot,  a meter, or other units, that fit into the region.  For the area between $y=f'(x)$ and $y=0$ over $[a, b]$, 
the simplest measure is to use an equivalent rectangle of length 
$L=b-a$ and width $W$ (See Figure 9). That is, the excess area above $y=W$ (the vertically stripped region) is moved to 
fill the deficit area (the horizontally striped region).  The corresponding description by a  mathematical formula  is below: 
\be
I[f'(x), a, b] = L\times W = (b-a)\times W.
\ee
This can be true as long as we have
\be
W= \frac{I[f'(x), a,b]}{b-a}=\frac{f(b)-f(a)}{b-a}.
\ee
We call this $W$  the mean value of $f'(x)$ over the interval $[a,b]$. 
Geometrically, $W=(f(b)-f(a))/(b-a)$ is the slope of the secant line that connects points $A$ and $B$ of Figure 10.
If $y=f(x)$ is not a straight line, then there must be a point $m$ in $[a,b]$ whose slope is less than $(f(b)-f(a))/(b-a)$,
and another point $M$ whose slope is larger than $(f(b)-f(a))/(b-a)$, i.e.,
\be
f'(m)\leq \frac{f(b-f(a)}{b-a}\leq f'(M).
\ee
Between $f'(m)$ and $f'(M)$, $(f(b)-f(a))/(b-a)$ must meet the mid-ground slope $f'(c)$ at a point $c$ in $[a,b]$, i.e.,
\be
W=f'(c) = \frac{f(b)-f(a)}{b-a}.
\ee
This is the mean value theorem (MVT) of calculus.   It states that 

\vskip 0.2in 
\noindent {\bf Theorem 1. (MVT).}
{\it There exists $c$ in $[a,b]$ such that $f'(c) = \frac{f(b)-f(a)}{b-a}$ if $f'(x)$ has a value for every $x$ in $[a,b]$.}
\vskip 0.2in

\begin{figure} [ht]
\centering
\includegraphics[height=3in,width=4in]{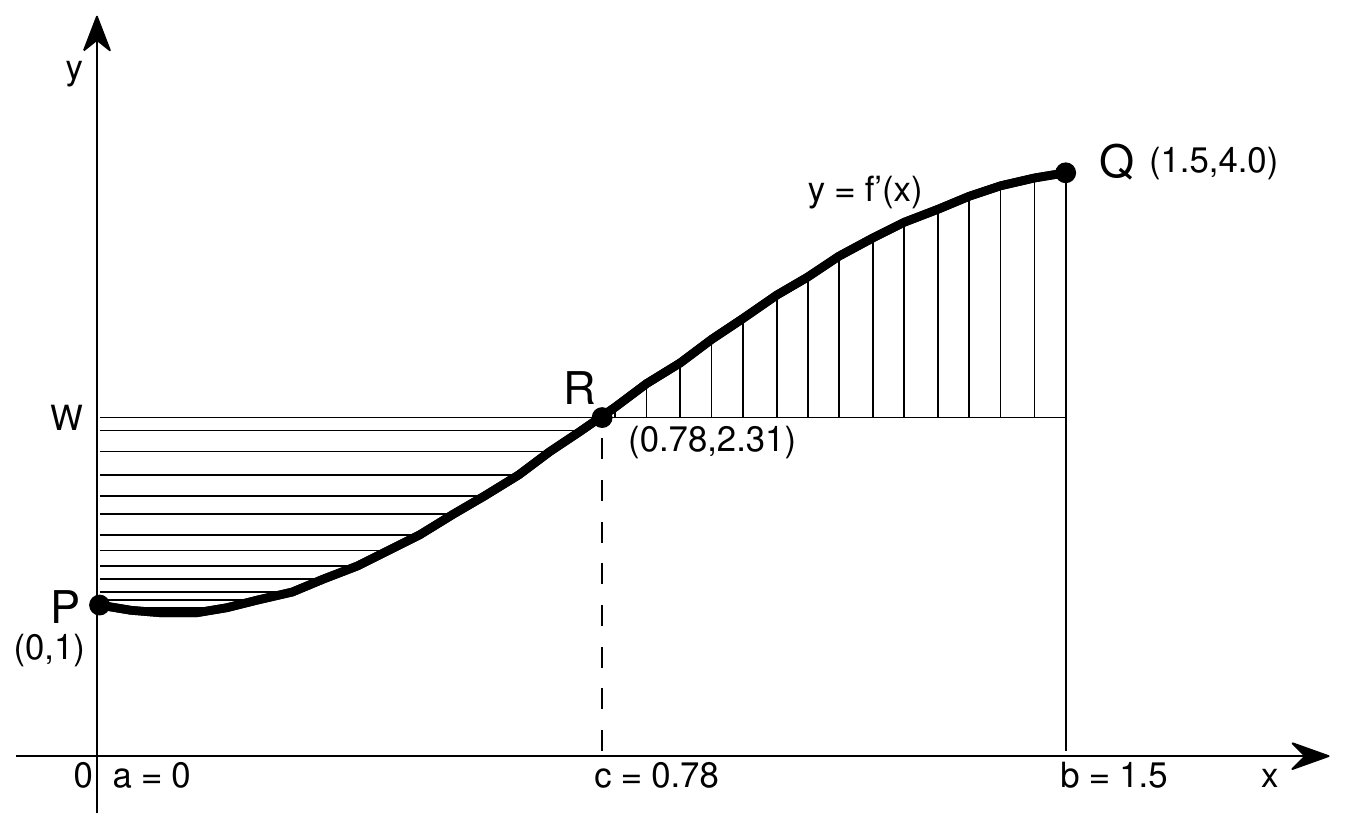} 
\caption{An area interpretation of an integral.  }
\label{figure9}
\end{figure}

\begin{figure} [ht]
\centering
\includegraphics[height=3in,width=4in]{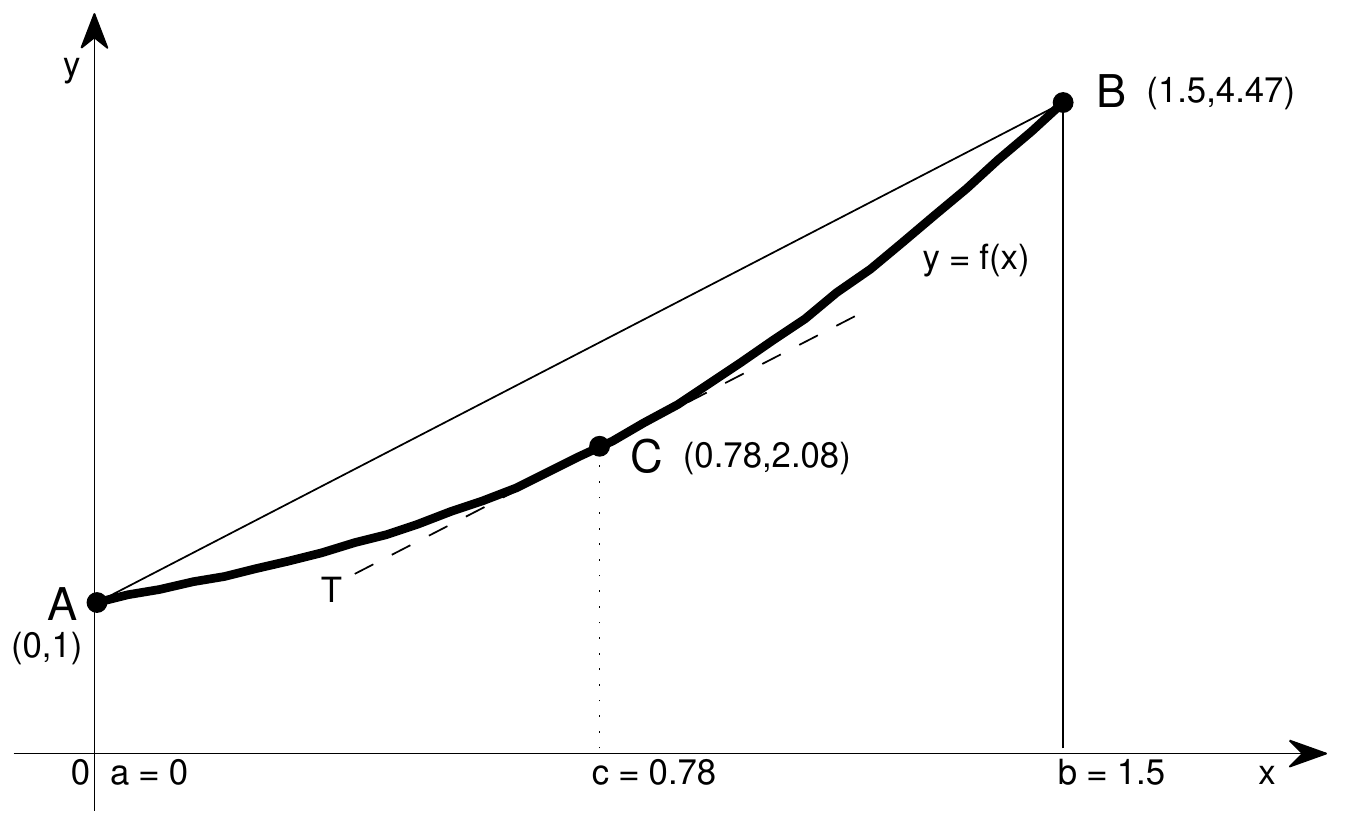} 
\caption{Illustration of the mean value theorem: There is a tangent line parallel to the secant line connecting points $A$ and $B$. }
\label{figure10}
\end{figure}

Geometrically, MVT means that there is at least one point $c$  whose tangent line is parallel to the secant line $AB$. Of course, this
holds if $y=f(x)$ is a straight line, in which case $c$ can be any point in $[a,b]$.

Rigorous mathematics for MVT would require one to prove the above statement ``$(f(b)-f(a))/(b-a)$ must meet the mid-ground slope $f'(c)$ at one point $c$ in $[a,b]$," 
 namely, it proves the existence of the point $c$. This is to prove the intermediate value theorem and is beyond the scope of this introductory lecture.

Therefore, the integral $I[f'(t),a,x] = f(x)-f(a)$ is the increment of the antiderivative from $a$ to $x$, and is also the area for the region between  
the integrand derivative function and $y=0$ in the interval $[a,x]$, i.e., the region bounded by $y=f'(t), y=0, t=a$ and $t=x$. 
 The traditional definition of an integral is from the aspect of an area that is defined as a sum of many rectangles of  increasing narrow widths, under the 
 condition of each width approaching zero. For non-mathematics majors and general public, the condition of each width approaching zero, which is  
 a concept of limit, adds complexity and confusion to the traditional definition of integral.   
 In contrast, the geometric meaning of our direct integral is the height increment of the antiderivative, not the 
 area underneath of the derivative. The area is only regarded as an additional geometric interpretation according to the 
 intermediate value theorem. Under this interpretation of area, we have the following example. 

\vskip 0.3cm

\noindent {\bf Example 7.} $I[\sqrt{1-x^2}, 0,1] $ is the area of a quarter unit round disc and is thus equal to $\pi/4$,
since $y=\sqrt{1-x^2}$ represents a quarter circle in the first quadrant. 

\vskip 0.3cm

Calculating the slope using the factorization method works for polynomial functions, but the procedure is tedious. The procedure may not even 
work for transcendental functions like $y=\sin x$.  The MVT provides another way to calculate the slope by using the slope of a secant line. 
In the above MVT, if $B$ moves very close to $A$, then the mean value $(f(b)-f(a))/(b-a)=f'(c)$ in Theorem 1 is approaching the slope at $A$, since $c$ is 
approaching $a$, forced by $b$ approaching $a$. The formal writing is 
\be
\lim_{b \to a}  \frac{f(b)-f(a)}{b-a} =f'(a). 
\ee
This  can also be considered  a definition of derivative and is   called {\it defining a derivative by a limit}. This procedure 
is efficient to calculate the derivative by hand and to derive many traditional derivative formulas in the earlier years of calculus.

\section{ Descartes' method of tangents and brief historical note}

\noindent  So, what is the origin of the direct calculus ideas described above? 
Numerous papers and books have discussed the historical development of calculus. Here we  recount
the development by  a few major historical figures  to sort out the origin of the main ideas of  the calculus 
method outlined above. 
We cite only a few modern references and two works by Newton. Our focus is on (i) Descartes' method of tangents
that is the earliest systematic way of finding the slope of a curve without using a limit, and (ii) Wallis' formulas 
of area, which were the earliest form of FTC. We do not intend to  present a complete 
list of important works on the history of calculus. 

In 1638, Ren{\' e} Descartes (1596-1650) derived his  method of tangents and 
 included the method in his 1649 book,  {\it Geometry} (see Cajori (1985), p.176). 
 Descartes' method of tangents is purely geometric, 
 constructing a tangent circle at a given point of a curve with the circle's center on
 the  x-axis (Cajori (1985), pp. 176-177). Graphically, it is easier to draw a tangent circle than a tangent 
  line using a compass and rulers. The tangent circle can  be constructed  using a radius 
line and moving the center on the x-axis so that the circle touches the curve at 
only one point. Then a tangent line can be drawn as the line perpendicular to the radial 
line of the circle at the tangent point. 

Descartes' method of tangents also has an analytic description. The tangent circle is determined by the given 
point $P(x_0, y_0)$ on the curve and the moving center on the x-axis $(a,0)$. The 
circle's equation is 
\be
\label{eq28}
(x-b)^2 + (y-0)^2 = (x_0-a)^2 + (y_0-0)^2.
\ee
The tangent condition requires that  this equation and the curve's equation $y=f(x)$ have a double root at $P(x_0, y_0)$. 
This can determine $a$ and hence the tangent circle. The radial line is determined by 
$P(x_0, y_0)$ and $(a, 0)$ and has its slope $m_R= y_0/(x_0-a)$. 
The tangent line of the circle at  point  $P(x_0, y_0)$ is the tangent line of the curve at the same point. 
The slope of the tangent 
line is calculated as $m_T=-1/m_R= (a-x_0)/y_0$.  In the above procedure, limit is not used. 

\vskip 0.3cm

\noindent {\bf Example 8. } Use Descartes' method of tangents to find the slope of $y=\sqrt{x}$ at $(1,1)$. 

Substituting $y=\sqrt{x}$ into eq. (\ref{eq28}), we obtain
\be
(x-a)^2 +x = (1-a)^2 + 1.
\ee
This can be simplified to
\be
x^2 + (1-2a)x + 2(a-1)=0.
\ee
Knowing that $x=1$ is a solution of this equation helps factorize the left-hand side
\be
(x-1) (x+2-2a)=0.
\ee
The double root condition for a tangent line requires $x_1=x_2=1$. This leads to 
\be
x_2+2-2a=1+2-2a=0.
\ee
Hence, $a=3/2.$
The slope of the radial line is $m_R= 1/(1-3/2)=-2$. The slope of the tangent line is thus $m_T=1/2$. 

\vskip 0.3cm

Although Descartes' method of tangents is complicated in calculation, its concept is 
simple, clear and unambiguous, and its geometric procedure is sound. It does  not involve  small increments of 
an independent variable (developed by Fermat also in the 1630s), and hence it does not involve limits or infinitesimals.  
According to the point-slope equation of a tangent line presented earlier, the complexity of Descartes' method of tangents
is unnecessary to calculate a slope. However, the 
point-slope form of a line was not known during Descartes' lifetime. According to Range [{\bf 6}], the  point-slope 
form of a line was first introduced explicitly by Gaspard Monge (1746-1818) 
in a paper published in 1784. Thus, Monge's point-slope method of tangent appeared more than 100 years after 
Descartes' method of tangents. 

Pierre de Fermat (1601-1665)'s method of tangents is similar to the modern method of differential quotient and uses a small 
increment (i.e., infinitesimal), which is ultimately set to be zero when the infinitesimal is forced to disappear from the denominator (Ginsburg et al., 1998).
Thus, Fermat's method of tangents is more efficient for calculation from the point of view of  limit, while Descartes'
method of tangents is geometrically more direct and easier to plot by hand, and Monge's method of 
tangents is geometrically more direct.
Fermat also used a sequence of rectangular strips to calculate the area 
under a parabola. His strips have variable width, which enabled him to use the sum of geometric series.
This method of calculating an area can be traced back to Archimedes. 

Archimedes (287-212 BC) 's method of exhaustion enabled him to find the area under a parabola. He used infinitely 
many triangles inscribed inside the parabola and also utilized the sum of   geometric series. 

Bonaventura Cavalieri (1589-1647) used  rectangular strips of equal width to calculate the area under a 
straight line (i.e., a triangle) and under  a parabola. 

Around 1655-1656, John Wallis (1616-1703) derived algebraic formulas that represent the areas under the curve of 
simple functions, such as $y=kt$ and $y=kt^2$, from $0$ to $x$ (Ginsburg et al., 1998). Considering the existing work on 
tangents (i.e., slopes or derivatives) at that time, and considering the DA pair concept here, 
 we thus may conclude that Wallis  had already explicitly demonstrated, before Newton, the relationship between slope and 
area using examples, i.e., FTC. 

Isaac Newton (1642-1727) attended Trinity College, Cambridge in 1660 and quickly made himself  a master of  Descartes' {\it Geometry}.
He learned much mathematics from his teacher and friend Isaac Barrow (1630-1677), who knew the method of tangents by both
Descartes and Fermat and also knew how to calculate areas under some simple functions. Barrow revealed 
that differentiation and integration were inverse operations, i.e., FTC (Cajori, 1985). Newton summarized the past work on tangents and area calculation,
introduced many applications of the two operations, and made the tangent and area methods a systematic set of mathematics theory. 
Newton's method of tangents followed that of Fermat and had a small increment that  eventually approaches zero. Namely, he used 
a sequence of secant lines to approach a tangent line as is done in  modern calculus' definition of derivative.  Although the ``method of limits" is
 frequently attributed to Newton (see his book  (Newton, 1729) entitled ``{\it The Mathematical 
Principles of Natural Philosophy}" (p45)), he was not as adamant as 
Leibniz about letting an infinitesimal be zero at the end of a calculation.  Newton was dissatisfied with the omitted small errors. 
He wrote that ``in mathematics the minutest errors are not to be 
neglected" (see Cajori (1985), p198). 

Newton's method of fluxions intended to solve two fundamental mechanics problems  that 
are equivalent to the two geometric problems of slope and height increment of a curve pointed out earlier in in  this paper:  
\\
``(i). The length of the space described being continually (i.e., at all times) given; to find the
    velocity of the motion at any time proposed.\\
   (ii). The velocity of motion being continuously given; to  find the length of the space described at any time proposed"
(see Cajori (1985), p.193, and  Newton (1736), p.19).
\\

The solution to these two problems also led to FTC, as geometrically interpreted in previous sections. 
The above statement of the two problems is directly cited from Newton's book  ``{\it The Method of Fluxions} "(Newton (1736), p.19),  
which was translated to English from Latin and published by John Colson. 

Gottfried Wilhelm Leibniz (1646-1716) produced a profound work similar to Newton's that  summarized the method of  tangents
and the method of area using a systematic approach. His approach has been passed on to today's classrooms, including his notations of derivative and 
integration. 

Our description of calculus method  has demonstrated that  if we avoid calculating the area underneath  a curve 
and define an integral by the height increment, we
can readily extend Descartes' method of tangents to establish the theory of differentiation and integration 
by considering the slope  (i.e., grade), DA pair, and  height increment.
The no-limit approach to  calculus  outlined in this paper is attributable to Descartes' original ideas, and 
is different from the  those of Fermat, Newton and Leibniz. The FTC is attributable to Wallis' original ideas. 
Ginsburg et al. (1998)  concluded that the query whether Leibniz 
plagiarized Newton's work on calculus is not 
really a valid question since the calculus ideas had already been developed by others 
before the calculus works of either Newton or Leibniz. Both just summarized the work of earlier mathematicians and 
developed  differentiation and integration into a systematic branch of mathematics by using
the methods of infinitesimals and limits. After their work, calculus became a very useful tool in engineering, 
natural sciences, and numerous other fields.

In addition to  the aforementioned mathematicians, there are many others who contributed to the development of calculus, including 
Gregory of St. Vincent (1584-1667) from Spain, Gilles Persone de Roberval (1602-1675), Blaise Pascal (1623-1662),
Christiaan Huygens (1629-1695), and Leonhard Euler (1707-1783). Augustin-Louis Cauchy (1789-1857) has been credited with the 
rigorous development of calculus from the definition of limits. Karl Weierstrass (1815-1897) corrected 
Cauchy's mistakes and introduced the delta-epsilon language we use today in mathematical 
analysis. 

\section{ Conclusions}

\noindent We have introduced the concepts of derivatives and integrals without using limits. Geometrically, derivatives were introduced 
directly from the slope of a tangent line. Algebraically, derivatives and antiderivatives were introduced simultaneously as a 
DA pair. Then the integral was 
introduced as the height increment of the antiderivative function. This increment was geometrically interpreted  as the area of the region 
bounded by the integrand function,
the horizontal axis and the integration interval. A justification of this interpretation was given to demonstrate that this definition of area was 
reasonable and mathematically rigorous up to the proof of IVT which, like the Euclidean axioms,  is intuitively true to most people 
in general public. At the end, we pointed out that limit was an efficient approach to calculate
derivatives by hand and could help derive derivative formulas for complicated functions besides polynomials. We thus regard
that the limit approach to calculus is an excellent computing method for finding derivatives by hand. In the pre-computer 
era, this limit approach was obviously critical in calculating derivatives of a variety of functions. In our current computer era, 
the limit approach to calculus is less essential and may be unnecessary for non-mathematical majors 
or general public at the introductory stage. 

Although  the ideas of direct calculus described in this paper come from practical applications, we have maintained the self-contained mathematical rigor and logic. 
Pure mathematical analysis regarding the structure of  real line and the sequence approach to a compact set are not topics
for this introduction. These analysis approaches mainly due to Cauchy and Weierstrass  have certainly enriched calculus as begun by
Archimedes, Descartes, Fermat, Wallis, Newton, Leibniz and others. However, our paper has shown that it is possible to introduce the basic concepts 
and calculation methods of calculus directly, without using limits. 

One can also regard calculus as an extension of the  trigonometry of a regular right triangle to the trigonometry of 
a curved-hypotenuse right triangle.  For a regular right triangle, the slope (i.e., derivative = tangent of the angle) of the hypotenuse  is 
derivative, and the vertical
increment (i.e., an  integral)   is equal to the opposite side, which is the integral of the derivative (see Figure 11). It is obvious that the vertical increment and the slope are related and 
have the following relationships:
\begin{equation}
\tan \theta = \frac{BC}{AC}  \qquad   \mbox{(derivative)}
\ee
and 
\begin{equation}
BC = \tan \theta \times AC   \qquad   \mbox{(integral)}. 
\ee
These two formulas  are the FTC for the regular right triangle. 
The extension is from this straight line hypotenuse to the curved hypotenuse,
such as a parabola or an exponential function. For the curved hypotenuse, the slope varies at different points and the 
total height  increment  is an integral. That is,
\be
m=f'(x) \qquad   \mbox{(derivative)}
\ee
and
\be
BC=I[f'(x), a,b]  \qquad   \mbox{(integral)}. 
\ee

\begin{figure} [ht]
\centering
\includegraphics[height=2.3in,width=5in]{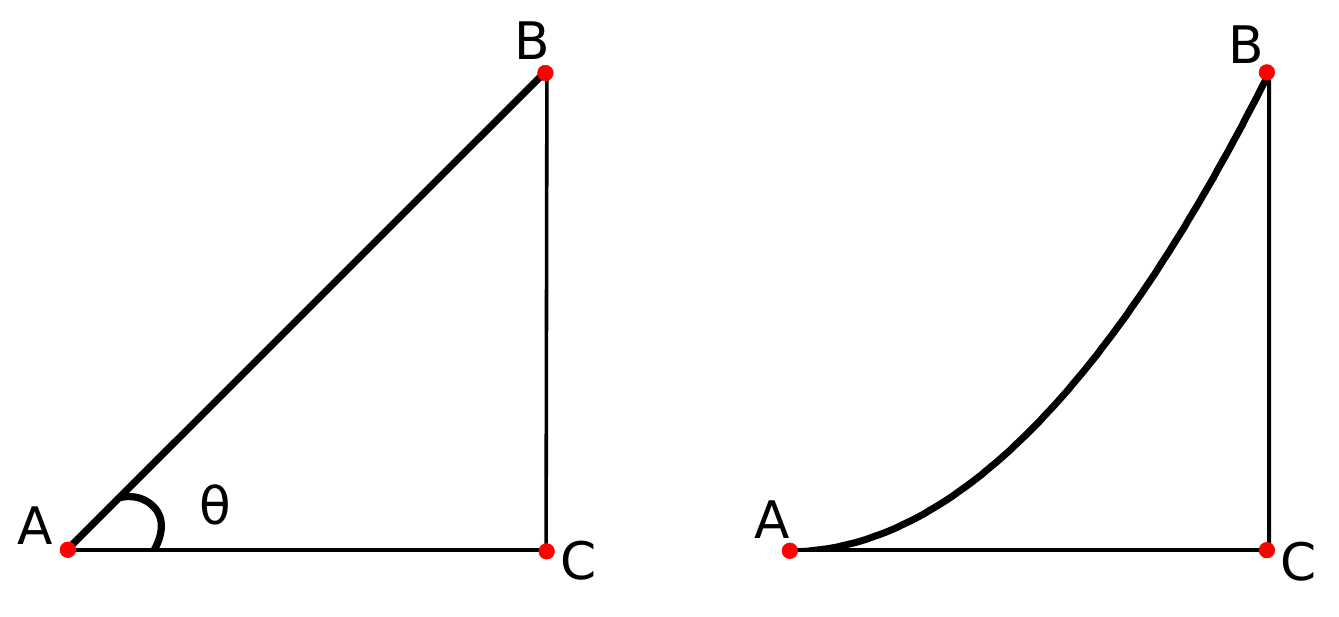}
\caption{Calculus explained by triangles with straight and curved hypotenuses }
\label{figure11}
\end{figure}

\vskip 0.5cm

\noindent  {\bf Acknowledgments.} ~ 
 Julien Pierret, Jazlynn Ngo, and Kimberly Leung assisted with plotting the figures in this paper. The discussion with 
 Chris Rasmussen and Dov Zakis was helpful in the paper's development. The work was partially supported by the 
 US National Science Foundation. During his many visits to China, the  lead author, Samuel S.P. Shen, exchanged ideas and worked together 
 with the second author, Qun Lin. The original ideas and technical details of this paper are credited to the lead author.

\vskip 1.0cm

\noindent {\Large\bf References}

\begin{enumerate}

\item F. Cajori, {\it A History of Mathematics } (pp. 162-198), 4th ed., Chelsea Publishing Co., New York, 534pp, 1985.

\item  J.L. Coolidge, The story of tangents. {\it American Math. Monthly} {\bf 58} (1951) 449-462. 

\item D. Ginsburg, B. Groose, J. Taylor, and B. Vernescu,  The History of the Calculus and the Development of Computer Algebra Systems, {\it Worcester Polytechnic 
Institute Junior-Year Project}, www.math.wpi.edu/IQP/BVCalcHist/calctoc.html, 1998 

\item Q. Lin,  {\it Calculus for High School Students: from a Perspective of Height Increment of a Curve}, People's Education Press, Beijing, 2010. 

\item Q. Lin, {\it Fastfood Calculus}, Science Press, Beijing, 2009. 

\item I. Newton, {\it The Method of Fluxions and Infinite Series with Application 
to the Geometry of Curve-Lines}, Translated from Latin to English by J. Colson, Printed for Henry Woodfall, London, 339pp, 1736. [http://books.google.com]
\begin{verbatim}
http://books.google.com/books?id=WyQOAAAAQAAJ&printsec=frontcover&source=
gbs_ge_summary_r&cad=0#v=onepage&q&f=false
\end{verbatim}

\item I. Newton, {\it The Mathematical 
Principles of Natural Philosophy}, Translated from Latin to English by A. Motte, Printed for Benjamin Motte, London, 320pp, 1729.  [http://books.google.com]
\begin{verbatim}
http://books.google.com/books?id=Tm0FAAAAQAAJ&printsec=frontcover&source=
gbs_ge_summary_r&cad=0#v=onepage&q&f=true
\end{verbatim}

\item R.M. Range, Where are limits needed in calculus? {\it Amer Math Monthly} {\bf 118} (2011) 404-417.

 \item J. Susuki, The lost calculus (1637-1670): Tangency and optimization without limits. {\it Mathematics Mag.} {\bf 78} (2005) 339-353.
\end{enumerate}

\newpage

\begin{figure} [ht]
\centering
\includegraphics[height=9.8in,width=7in]{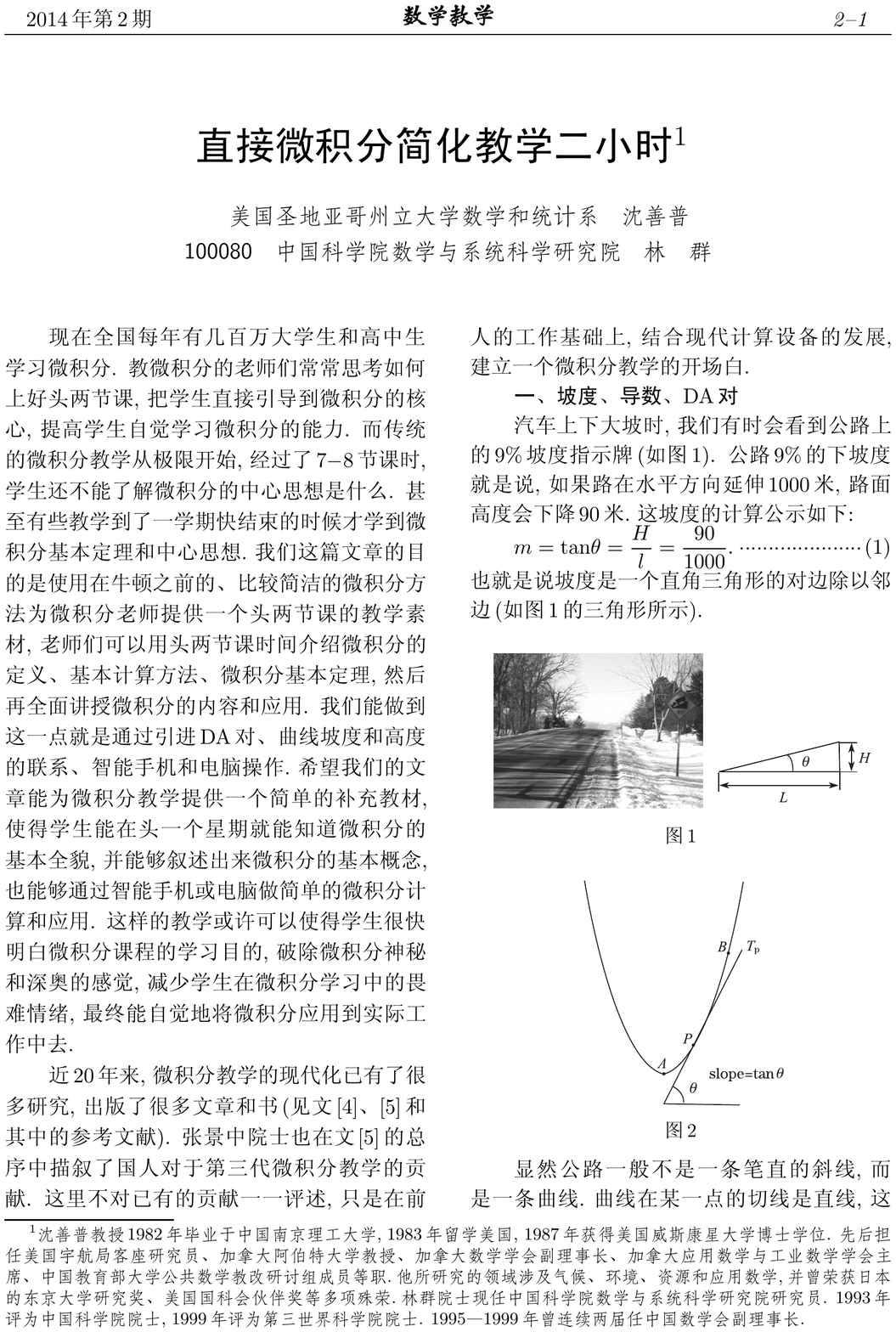}
\label{figure11}
\end{figure}

\end{document}